

\baselineskip=14pt
\parskip=10pt
\def\halmos{\hbox{\vrule height0.15cm width0.01cm\vbox{\hrule height
  0.01cm width0.2cm \vskip0.15cm \hrule height 0.01cm width0.2cm}\vrule
  height0.15cm width 0.01cm}}

\magnification=\magstephalf

\def\1{{\overline{1}}}
\def\2{{\overline{2}}}
\parindent=0pt
\overfullrule=0in

\def\frac#1#2{{#1 \over #2}}
\centerline
{\bf  Bijective and Automated Approaches to Abel Sums }
\bigskip
\centerline
{\it Gil KALAI and Doron ZEILBERGER}

\qquad \qquad {\it  Dedicated to Dominique Foata (b. Oct. 12, 1934), on his forthcoming 90th birthday}

{\bf Abstract}: In this tribute to our guru, Dominique Foata, we return to one of our (and Foata's) first loves, and revisit 
Abel sums and their identities, from two different viewpoints.

{\bf Preface}

In the very first issue of {\it Crelle's journal} (the first mathematical periodical solely dedicated to mathematics), Niels Henrik Abel published
a two-page paper [A], stating and proving his eponymous identity. This lead to an intensive study of general {\it Abel Sums} by many people, 
(see [R] [C] and their numerous references), and
to beautiful bijective approaches pioneered by Dominique Foata and Aim\'e Fuchs [Fo][FoFu] that led to Francon's elegant proof [Fr] (see [C], p. 129).
This tribute consists of two independent parts. The first part is {\it bijective}, while the second
part is {\it automated}, elaborating and extending  John Majewicz's 1997 PhD thesis [M1][M2], and more important,
fully implementing it (and its extension) in a Maple package   \hfill\break
{\tt https://sites.math.rutgers.edu/\~{}zeilberg/tokhniot/AbelCeline.txt} \quad .

{\bf Part I: Bijective Proofs (\`a la Foata) of an Abel-Type identity and a generalization}

In [Ka] (see also [Ka']), the first author proved (as a special case of  more general results) the following Abel-type identity. Let $n,p$ by nonnegative integers, then

$$
\sum _{k=0}^n {{n} \choose {k}} k^k (n-k)^{n-k+p} = \sum _{k=0}^n {{n} \choose {k}}n^k (n-k)! S(p+n-k,n-k) \quad .
\eqno(1)
$$

Here, $S(n,k)$ are the Stirling numbers of the second kind, and $k! S(n,k)$ is the number of maps from a set of size $k$ onto a set of size $n$. (This property can serve as the definition of these numbers.) 

A special case of this formula is attributed to Cauchy (Equation (24) in Chapter 1 of [R]):

$$
\sum _{k=0}^n {{n} \choose {k}} k^k (n-k)^{n-k} = \sum _{k=0}^n {{n} \choose {k}}n^k (n-k)! \quad .
\eqno(2)
$$

Here we present a combinatorial proof, in the style of Foata [Fo][FoFu] for Formulas $(2)$ and $(1)$. 
We use the notation $[n]=\{1,2,\dots,n\}$, and to make the argument clearer we will present both proofs.

\vfill\eject

{\bf The Proofs}

{\bf Proof of Equation $(2)$}:
The left hand side counts triples $f,A,B$ where
$f$ is a function from [n] to [n]
and the following conditions hold:

$$
A \cup B=[n], A \cap B =\emptyset, f(A)\subset A {\rm ~and~}  f(B) \subset B.
\eqno(3)
$$

Indeed, if $|A|=k$, there are ${{n} \choose {k}}$
ways to choose $A$ ($B$ is now determined, $B=[n] \backslash A$), $k^k$ ways to choose the restriction of $f$ to $A$ subject to the condition $f(A) \subset A$ and $(n-k)^{n-k}$ ways to choose the restriction of $f$ to $B$ subject to the condition $f(B) \subset B$.  

The right hand side counts triples $f,C,D$ where
$f$ is a function from $[n]$ to $[n]$
and following conditions hold:

$$
C \cup D=[n], C \cap D =\emptyset, f(D) = D.
\eqno(4)
$$
(Here, 
$f(D)=D$ means that $f$ is a bijection from $D$ to $D$; Note that we relaxed the condition for $C$ compared to $A$ and strengthened the condition for $D$ compared to $B$.) Indeed, if $|C|=k$, there are ${{n} \choose {k}}$
ways to choose $C$, $n^k$ ways to choose the restriction of $f$ to $C$ (no conditions here) and  $(n-k)!$ ways to choose the restriction of $f$ to $D$ subject to the condition $f(D) =D$.

The crucial observation is:

$\bullet$ For {\it every} function $f$ from $[n]$ to $[n]$ the number of pairs $(A,B)$ that satisfies Equation $(3)$ equals the number of pairs $(C,D)$ that satisfies equation $(4)$. 

Indeed if $(A,B)$ is a pair that satisfies Equation $(3)$, we can take 
$$
D=f(f(\cdots (f(B))\cdots ) {\rm ~and~} C=f^{-1}(f^{-1} (\cdots (f^{-1} (A))\cdot ).
\eqno(5)
$$

In other words $D=\{d \in [n]: d=f^k(a) {\rm ~for~some~} a \in B, {\rm~ and~} k \ge 0\}$ and $C=\{c \in [n]: f^k(c)\in B {\rm ~for~some~} k \ge 0\}$. Taking the inverse operation brings you from a pair $(C,D)$ to a pair $(A,B)$

{\bf Proof of Equation $(1)$}:

Let $X=[n]=\{1,2,\dots. n\}$, and $Y=[n+p]$ 
Consider all triples $(f, A, B)$ where, this time, $f$ maps $[n+p]$ to $[n]$, and $A,B$ satisfy

$$
A \cup B=[n], A \cap B =\emptyset, f(A)\subset A {\rm ~and~}  f(B \cup [n+1,n+p]) \subset B.
$$

We also consider triples $(f, C, D)$ where, $f$ maps $[n+p]$ to $[n]$, and $C,D$ satisfy

$$ 
C \cup D=[n], C \cap D =\emptyset, f(D \cup [n+1,n+p]) = D.
\eqno(7)
$$

Also in this case a stronger statement holds: For every function  $f: [n+p] \to [n]$  there is a bijection between between pairs $(A,B)$ satisfying equation $(6)$ and pairs $(C,D)$ satisfying equation 
$(7)$ and the bijection is given again by $(5)$.

This bijection is implemented in the Maple package {\tt AbelBijection.txt}, available from

{\tt https://sites.math.rutgers.edu/\~{}zeilberg/tokhniot/AbelBijection.txt} \quad .

See the front of this article \hfill\break
{\tt https://sites.math.rutgers.edu/\~{}zeilberg/mamarim/mamarimhtml/abelKZ.html} for a sample input file and its corresponding output file.

{\bf Remark:} If we give each element $k \in [n]$ a weight $w_k$ and a function $f$ from $[n+p] \to [n]$ the weight $w(f):=\prod \{w_{f(k)} \, : \, k \in [n]\}$ then our bijective proof 
gives a Hurwitz-type generalization of Abel's formula, see Hurwitz [H] and Excercise 30 in Section 2.3.4.4. of Knuth [Kn]. 
(For Formula $(1)$ we obtain a Hurwitz-type generalization of the Stirling numbers of the second kind.) 
 
{\bf Part II: Automating Abel Sums}

Abel's original identity had many proofs, but the one by Shalosh B. Ekhad and John Majewicz [EM] (written $30$ years ago, and dedicated to Dominique Foata on his $60^{th}$ birthday)
has the distinction that it was {\it computer-generated}, yet with a bit of patience,
{\it humanly readable} and verifiable. In order to motivate the sequel, let's reproduce it in its entirety.

{\bf Abel's identity}: For any non-negative integer, $n$
$$
\sum_{k=0}^n {{n} \choose {k}} \,(r+k)^{k-1}(s-k)^{n-k}=\frac{(r+s)^n}{r}  \quad .
\eqno(8)
$$
 
{\bf Proof}  Let $F_{n,k}(r,s)$ and $a_n(r,s)$ denote, respectively, 
the summand and sum on the LHS of (8), and let 
$G_{n,k}:=(s-n) {{n-1} \choose {k-1}} \,(k+r)^{k-1}(s-k)^{n-k-1}$. Since
$$
F_{n,k}(r,s)-sF_{n-1,k}(r,s)-(n+r)F_{n-1,k}(r+1,s-1)
        +(n-1)(r+s)F_{n-2,k}(r+1,s-1)=G_{n,k}-G_{n,k+1} \,,
$$
(check!), we have by summing from $k=0$ to $k=n$, thanks to the telescoping
on the right:
$$
 a_n(r,s)-sa_{n-1}(r,s)-(n+r) a_{n-1}(r+1,s-1)+
        (n-1)(r+s) a_{n-2}(r+1,s-1)=0.
\eqno(9)
$$
Since $(r+s)^n\cdot r^{-1}$ also satisfies 
this recurrence (check!) with the same initial 
conditions $a_0(r,s)=r^{-1}$ and 
$a_1(r,s)=(r+s)\cdot r^{-1}$, $(8)$ follows.\halmos

This proof was derived using John Majewicz's brilliant adaptation of {\it Sister Celine's technique} [Fa1][Fa2][Z] (see also [PWZ], chapter 4). Recall that Sister Celine was interested in finding {\it pure recurrence relations} of the form
$$
c_0(n)a_n+  c_1(n)a_{n+1}+  \dots + c_L(n) a_{n+L}=0 \quad,
\eqno(10)
$$
where $c_0(n), \dots, c_L(n)$ are polynomials in $n$ for sequences $a_n$, that are  defined by expressions of the form
$$
a_n:=\sum_{k=-\infty}^{\infty} F_{n,k} \quad,
$$
where  $F_{n,k}$ is proper hypergeometric (see [PWZ] for the definition, in particular,  $F_{n+1,k}/F_{n,k}$ and $F_{n,k+1}/F_{n,k}$ are both rational functions of $n$ and $k$) .
The way she did it was to (in her case by hand, but nowadays it has all been fully automated) to search for a recurrence of the form: 
$$
\sum_{i=0}^{L} \sum_{j=0}^{M} b_{ij}(n) F_{n+i,k+j} \, = \, 0 \quad ,
\eqno(11)
$$
for  {\it some} positive integers $L$ and $M$.

Dividing by $F_{n,k}$, and clearing denominators, looking at the numerator,  and then setting all the coefficients of powers of $k$ to $0$, yields a system of {\it linear equations} (with coefficients that
are polynomials in $n$) for the {\it undetermined} $b_{ij}(n)$.

Having found the $b_{ij}(n)$, summing Eq. $(11)$ from $k=-\infty$ to $k=\infty$, yields Eq. $(10)$ with
$$
c_i(n)=\sum_{j=0}^{M} b_{ij}(n) \quad, \quad 0 \leq i \leq L \quad .
$$

In his  PhD thesis  [M1][M2] (written under the direction of the second author), John Majewicz adapted Sister Celine's method to {\it Abel Type} sums, of the form
$$
a_n(r,s)=\sum_{k=0}^n \,F_{n,k} \, (r+k)^{k-1+p}(s-k)^{n-k+q} x^k  \quad ,
\eqno(12)
$$
where $F_{n,k}$ is hypergeometric in $n$ and $k$. Here $p$ and $q$ are arbitrary integers, and $x$ is any number (or symbol). It is no longer the case that
$a_n(r,s)$ satisfies a {\it pure} recurrence in $n$, with $r$ and $s$ {\bf fixed}, but it does satisfy a {\it functional recurrence}, similar (but often much more complicated) to Eq. $(10)$.
Denoting the summand of $(12)$ by $\overline{F}_{n,k}(r,s)$
$$
\overline{F}_{n,k}(r,s) \, := \, F_{n,k} (r+k)^{k-1+p}(s-k)^{n-k+q} x^k  \quad ,
$$
one looks for polynomials $b_{ij}(n)$ (that also depend on $r,s,p,q$ and $x$, but are {\bf free} of $k$), such that
$$
\sum_{i=0}^{L} \sum_{j=0}^{M} b_{ij}(n) \overline{F}_{n+i,k+j}(r-j,s+j) \, = \, 0 \quad .
\eqno(13)
$$
Dividing by  $\overline{F}_{n,k}(r,s)$ (since $\overline{F}_{n+i,k+j}(r-j,s+j)/\overline{F}_{n,k}(r,s)$ is still a rational function of $n$ and $k$), clearing denominators, looking at the
numerator, and setting all the coefficients of powers of $k$ to $0$, we get again a system of linear equations for the undetermined quantities $b_{ij}(n)$. We then ask Maple
to kindly {\tt solve} them, and if {\it in luck} we get a non-zero solution. It can be shown ([M1][M2]) that one can always find orders $L$ and $M$ for which such a system
is solvable (for sufficiently large $L$ and $M$, there are more unknowns than equations). Having found such a recurrence for $\overline{F}_{n,k}$, summing over $k$ yields
a functional recurrence for $a_n(r,s)$.
$$
\sum_{i=0}^{L} \sum_{j=0}^{M} b_{ij}(n) a_{n+i}(r-j,s+j) \, = \, 0 \quad .
\eqno(14)
$$

{\bf Implementation}

We fully implemented the Celine-Majewicz algorithm in a Maple packge {\tt AbelCeline.txt}  available from

{\tt https://sites.math.rutgers.edu/\~{}zeilberg/tokhniot/AbelCeline.txt} \quad . The function call is

{\tt FindOpe(F,n,k,r,s,R,S,N,K,MaxOrd); } \quad,

where {\tt MaxOrd} is the maximal order of the recurrence you are willing to tolerate.

It inputs a hypergeometric term $F_{n,k}$ and outputs the recurrence in the form of the corresponding operator, where $N$, $R$, and $S$ are the forward shift operators in
$n$, $r$, and $s$, respectively. To get a computer-generated paper, in {\it humanese}, the function call is:

{\tt Paper(F,n,k,r,s,R,S,N,K,MaxOrd): } \quad .

{\bf Sample Output}

Typing

{\tt Paper(binomial(n,k)*x**k,n,k,r,s,2):},

yields in $0.12$ seconds the following deep fact.

Let, for any integers $p$ and $q$ and number (or symbol) $x$
$$
a_n(r,s) := \sum_{k=0}^{n} {{n} \choose {k}} (r+k)^{k-1+p}(s-k)^{n-k+q}\,x^k \quad,
$$
then
$$
a_{n}\! \left(r , s\right) = 
\left(n x +r x \right) a_{n -1}\! \left(r +1, s -1\right)+s a_{n -1}\! \left(r , s\right)+\left(-n r x -n s x +r x +s x \right) a_{n -2}\! \left(r +1, s -1\right) \quad .
$$

Typing

{\tt Paper(1/(k!**2*(n-k)!)*x**k,n,k,r,s,3):}

yields that the sequence of polynomials in $(r,s)$ defined by
$$
a_n(r,s) := \sum_{k=0}^{n} \frac{1}{k!^2(n-k)!} \, (r+k)^{k-1+p}(s-k)^{n-k+q}\,x^k \quad,
$$
satisfies the functional recurrence
$$
a_{n}\! \left(r , s\right) = 
\frac{x \left(n +r \right) a_{n -1}\! \left(r +1, s -1\right)}{n^{2}}+\frac{\left(2 n^{2}-2 n s -2 n +s \right) s a_{n -1}\! \left(r , s\right)}{\left(n -1-s \right) n^{2}}
$$
$$
-\frac{\left(n^{2}+2 n r -2 r s -s^{2}-n -r \right) x a_{n -2}\! \left(r +1, s -1\right)}{\left(n -1-s \right) n^{2}}
$$
$$
-\frac{\left(n -s \right) s^{2} a_{n -2}\! \left(r , s\right)}{\left(n -1-s \right) n^{2}}+\frac{\left(n r +n s -r s -s^{2}\right) x a_{n -3}\! \left(r +1, s -1\right)}{\left(n -1-s \right) n^{2}} \quad .
$$

To see numerous other examples, read  the output files in the front of this article:

{\tt https://sites.math.rutgers.edu/\~{}zeilberg/mamarim/mamarimhtml/abelKZ.html} \quad .

In particular, to see the functional recurrence satisfied by the `innocent' sum:
$$
a_n(r,s) := \sum_{k=0}^{n} {{n} \choose {k}}  {{n+k} \choose {k}} (r+k)^{k-1+p}(s-k)^{n-k+q}\,x^k \quad,
$$
see the output file

{\tt https://sites.math.rutgers.edu/\~{}zeilberg/tokhniot/oAbelCeline4.txt} \quad .

Notice how complicated things get!

{\bf Getting differential recurrences}

Since $r$ and $s$, as opposed to $n$ and $k$, are `continuous variables', and since the functional recurrences gotten by Majewicz' Abel-Celine technique get so
complicated even with very simple summands, it occurred to us to look for {\it differential} recurrences.  For any bi-variate (proper) hypergeometric term $F_{n,k}$, defining (as above)
$$
a_n(r,s) := \sum_{k=0}^{n} F_{n,k} (r+k)^{k-1+p}(s-k)^{n-k+q}\,x^k \quad,
$$
where $p$ and $q$ are integers (but may be left symbolic) and $x$ is a number (again, it can be left symbolic),
one looks for differential-recurrence equations of the form
$$
\sum_{i=0}^{L}\sum_{j=0}^{M} b_{i,j}(n,r,s) \frac{d^{i}}{dr^i} a_{n+j}(r,s) \, = \, 0 \quad ,
$$
and
$$
\sum_{i=0}^{L}\sum_{j=0}^{M} c_{i,j}(n,r,s) \frac{d^{i}}{ds^i} a_{n+j}(r,s) \, = \, 0 \quad ,
$$
Together these pairs of differential recurrence equations, combined with initial conditions, uniquely determine $a_n(r,s)$.

In fact, one can consider far more general `kernels'. We can find a pair of differential recurrence equations for sums of the from
$$
a_n(r,s) := \sum_{k=0}^{n} F_{n,k} K(r,s,n,k) \quad,
$$
for any {\it kernel} $K(r,s,n,k)$ such that $(\frac{\partial K}{\partial r})/K$, and  $(\frac{\partial K}{\partial s})/K$ are rational functions of $(r,s)$ (and $k$).

We proceed analogously. One applies the generic operator to the summand $\overline{F}(n,k)$, expands, clears denominators, looks at the numerator, then equates
all the powers of $k$ to $0$, getting a system of linear equations that Maple can {\tt solve} for you.

To get these pairs of equations, in verbose, human-readable from, type:

{\tt PaperD(F,n,k,r,s,MaxOrd,KER):}

where {\tt MaxOrd} is the  maximum order you are willing to tolerate, and {\tt KER} is the kernel. For example, typing

{\tt PaperD(binomial(n,k),n,k,r,s,2,(r+k)**(k-1+p)*(s-k)**(n-k+q)*x**k):} \quad,

yields in a fraction of a second the facts that the sequence of polynomials in $r$ and $s$, $a_n(r,s)$ (for any $p$ and $q$ and $x$)
$$
a_n(r,s) := \sum_{k=0}^{n} {{n} \choose {k}} (r+k)^{k-1+p}(s-k)^{n-k+q}\,x^k \quad,
$$
satisfy the pair of recurrence-differential equations (the first in $r$, the second in  $s$)
$$
-\left(p n +n s -n +p +s -1\right) a_{n}\! \left(r , s\right)+\left(n r +n s +r +s \right) \left(\frac{\partial}{\partial r}a_{n}\! \left(r , s\right)\right)
$$
$$
+(n+p)a_{n +1}\! \left(r , s\right)-\left(n +r +1\right) \left(\frac{\partial}{\partial r}a_{n +1}\! \left(r , s\right)\right)=0 \quad,
$$
and
$$
-\left(n +1\right) \left(q +n -s +1\right) a_{n}\! \left(r , s\right)+q a_{n +1}\! \left(r , s\right)+\left(n -s +1\right) \left(\frac{\partial}{\partial s}a_{n +1}\! \left(r , s\right)\right) \, = 0 \quad .
$$

Note that setting $x=1$ and $p=0,q=0$ yields yet-another (automatic) proof of the original identity $(8)$. In fact it is closer in spirit to Niels Abel's original proof, that also used
differentiation (or rather integration).

{\bf Remark:} The proof technique of this part can prove (2) and can be adapted (with some pre-processing) to prove (1).

To see many other examples, look at the output files

{\tt https://sites.math.rutgers.edu/\~{}zeilberg/tokhniot/oAbelCeline5.txt} \quad ,

and

{\tt https://sites.math.rutgers.edu/\~{}zeilberg/tokhniot/oAbelCeline6.txt} \quad .

To see examples with more complicated kernels, see

{\tt https://sites.math.rutgers.edu/\~{}zeilberg/tokhniot/oAbelCeline7.txt} \quad .

Readers can easily find many other such deep  facts by playing with {\tt AbelCeline.txt}. Enjoy!

{\bf References}

[A] Niels H. Abel, {\it Beweis eines Ausdruckes, von welchem die Binomial-Formel ein einzelner Fall ist},   Crelle's J. Mathematik {\bf 1} (1826), 159-160. \hfill\break
{\tt https://sites.math.rutgers.edu/\~{}zeilberg/akherim/Abel1826.pdf} \quad .

[C] Louis Comtet, {\it ``Advanced Combinatorics''}, D. Reidel Publ. Co.,  Dordrecht/Boston, 1974 \quad . 

[EM] Shalosh B. Ekhad and John Majewicz, {\it A short WZ-style proof of Abel's identity}, Electronic Journal of Combinatorics {\bf 3(2)} (the Foata  Festschrift) \#R16. \hfill\break
{\tt https://sites.math.rutgers.edu/\~{}zeilberg/mamarim/mamarimPDF/abel.pdf} \quad .

[Fa1] Sister Mary Celine Fasenmyer, {\it Some generalized hypergeometric polynomials}, Bull. Amer. Moth. Soc. {\bf 53} (1947), 806-812.

[Fa2]  Sister Mary Celine Fasenmyer,  {\it A note on pure recurrence relations}, Amer. Math. Monthly {\bf 56} (1949), 14-17.

[Fo]  Dominique Foata, {\it Enumerating $k$-Trees}, Discrete Mathematics {\bf 2}  (1971), 181-186.\hfill\break
{\tt https://irma.math.unistra.fr/\~{}foata/paper/pub16.pdf} \quad .

[FoFu]  Dominique Foata and Aim\'e Fuchs, {\it  R\'earrangements de fonctions et d\'enombrement}, Discrete Mathematics {\bf 2}  (1971), 181-186. \hfill\break
{\tt https://irma.math.unistra.fr/\~{}foata/paper/pub12.pdf} \quad .

[Fr] Jean Francon, {\it Preuves combinatoires des identit\'es d'Abel}, Discrete Mathematics  {\bf 8} (1974), 331-343.

[H] Adolf Hurwitz, {\it \"Uber Abel’s Verallgemeinerung der binomischen Formel,} Acta Math. 26 (1902), 199-203.

[Ka] Gil Kalai, {\it A note on an evaluation of Abel sums},  J. Combin, Theory (Ser. A), {\bf 27} (1979), 213-217.\hfill\break
{\tt https://sites.math.rutgers.edu/\~{}zeilberg/akherim/kalai1979.pdf} \quad .

[Ka'] Gil Kalai, Nostalgia corner: John Riordan's referee report of my first paper, blog post on 
Combinatorics and More, 2021. https://gilkalai.wordpress.com/2021/02/19/nostalgia-corner-john-riordans-referee-report-of-my-first-paper/ 

[Kn] Donald E. Knuth, {\it The Art of Computer Programming,} Vol. I: Fundamental Algorithms (Addison-Wesley, Reading, MA, 1969).

[M1]  John Majewicz, {\it WZ-type certification procedures and Sister Celine's technique for  Abel-type sums}, J. of Difference Equations and Applications {\bf 2} (1996),  55-65.

[M2]  John Majewicz, {\it WZ certification of Abel-type identities and Askey's positivity conjecture}, Ph.D. thesis, Temple University 1997. \hfill\break
{\tt https://sites.math.rutgers.edu/\~{}zeilberg/Theses/JohnMajewiczThesis.pdf} \quad .

[PWZ] Marko Petkovsek, Herbert S. Wilf, and Doron Zeilberger, {\it ``A=B''}, A.K. Peters, 1996. \hfill\break
{\tt https://sites.math.rutgers.edu/\~{}zeilberg/AeqB.pdf} \quad .

[R] John Riordan, {\it ``Combinatorial identities''}, Wiley, New York, 1968.

[Z] Doron Zeilberger,  {\it Sister Celine's technique and its generalizations}, 
J. Math. Anal. Appl. {\bf 85} (1982), 114-145. \hfill\break
{\tt https://sites.math.rutgers.edu/\~{}zeilberg/mamarimY/celine1982.pdf} \quad .

\bigskip
\hrule
\bigskip
Gil Kalai, Hebrew University of Jerusalem, Einstein Institute of Mathematics, and
Reichman University, Efi Arazi School of Computer Science. \hfill\break
Email: {\tt gil dot kalai at gmail dot com} \quad .

\bigskip
Doron Zeilberger, Department of Mathematics, Rutgers University (New Brunswick), Hill Center-Busch Campus, 110 Frelinghuysen
Rd., Piscataway, NJ 08854-8019, USA. \hfill\break
Email: {\tt DoronZeil at gmail  dot com}   \quad .
\bigskip
{\bf May 1, 2024} \quad .

\end